**Integer Solutions, Rational solutions of the equations**

$$x^4 + y^4 + z^4 - 2x^2y^2 - 2y^2z^2 - 2z^2x^2 = n$$
$$\text{and } x^2 + y^4 + z^4 - 2xy^2 - 2xz^2 - 2y^2z^2 = n;$$

**And Crux Mathematicorum Contest**

**Corner problem CC24**

*Konstantine Zelator*

*P.O. Box 4280*

*Pittsburgh, PA 15203*

*U.S.A.*

*Email address: konstantine_zelator@yahoo.com*

*Konstantine Zelator*

*Department of Mathematics*

*301 Thankeray Hall*

*University of Pittsburgh*

*Pittsburgh, PA 15260*

*U.S.A.*

*Email address: spaceman@pitt.edu*

### 1. Introduction

In the May 2013 issue of the journal Crux Mathematicorum (vol. 38, No 5, see reference [1]), the following problem appeared:
(Contest Corner) CC24

1. Show that the equation $x^4 + y^4 + z^4 - 2x^2y^2 - 2y^2z^2 - 2z^2x^2 = 24$ has no integer solutions.

2. Does it have any rational solutions? Find one, or show that the above equation has no solutions.

Contest Corner problem CC24 is the motivation behind this work. In this article, we study the equation, in three variables x,y,z;

$$x^4 + y^4 + z^4 - 2x^2y^2 - 2y^2z^2 - 2z^2x^2 = n, \quad (1)$$

where n is a fixed or given positive integer.

We provide a soluition to the first questions of CC24, in Theorem 1: we show that if n is of the form, n=8N, where N is an odd positive integer then equation (1) has no integer solutions. So question 1 above, is really a particular case of Theorem 1; the case N=3.
Using Theorem 2, we answer question 2 of CC24. Theorem 2 postulates that if the integer n satisfies certain divisor conditions; then rational solutions do exist. Based on Theorem 2, for n=24, we find the rational solution $(x, y, z) = (\frac{5}{2}, \frac{1}{2}, 1)$ Obviously, since equation (1) is symmetric with respect to the three variables and since all the exponents are even; one really obtains 48=(6)(8) rational solutions from the above solution:

$$\left(\pm\frac{5}{2}, \pm 1, \pm\frac{1}{2}\right), \left(\pm\frac{5}{2}, \pm\frac{1}{2}, \pm 1\right), \left(\pm 1, \pm\frac{5}{2}, \pm\frac{1}{2}\right);$$

$$\left(\pm 1, \pm\frac{1}{2}, \pm\frac{5}{2}\right), \left(\pm\frac{1}{2}, \pm 1, \pm\frac{5}{2}\right), \left(\pm\frac{1}{2}, \pm\frac{5}{2}, \pm 1\right).$$

With all the sign combinations being possible (see Remark 1 for more details).
In Theorem 3, we prove a necessary and sufficient criterion for equations (1) to have an integer solution. Based on Theorem 3, in Theorem 4 we prove that equation (1) has no integer solutions if n=p (a prime), n=4, or n=p·q; where p and q are distinct primes.

Theorem 5 gives, in essence, a method for finding the nonnegative integer solutions of equation (1).
In Theorem 6, we find all the integer solutions to equation (1), twenty four in total; in the case $n = p^2$, p an odd prime.

In Theorem 7, we determine the integer and the rational solutions of equation (1) in the case $n = k^2$, k a positive integer; and with two among x, y, z being zero.

Theorem 8 lays the case $n = k^2$, $k\varepsilon\ \mathbb{Z}^+; k \geq 2$. and with exactly one of the three variables x, y, z being zero. All such integer solutions are determined.

Similarly, Theorem 9 deals with the case $n = k^2$, $k \geq 2$ and with exactly one of x,y,z being zero. All such rational solutions are determined.

Theorem 10 deals with the case n=1. All integer and rational solutions are listed.

Theorems 11 through 13 deal with the equation (quadratic in x), $x^2 + y^4 + z^4 - 2xy^2 - 2xz^2 - 2y^2z^2 = n$
In Theorem 11, the integer solutions of this equation (equation (11)), are described in detail.

Theorem 12 states that if n=2 or 3(mod4); then the above (quadratic in x) equation, has no integer solutions. Finally, Theorem 13 parametrically describes all the rational solutions of equation (11).

## 2. Three Lemmas

**Lemma1**

*Over the integers, and also over the rationals; equation (1), and equations (2) and (3) below, have the same solution set. In other words equations (1), (2), and (3) are equivalent both over the rationals and over the integers.*

$$(x+y+z) \cdot (x+y-z) \cdot (x-y+z) \cdot (x-y-z) = n \quad (3)$$

$$(x^2 + y^2 + 2xy - z^2) \cdot (x^2 + y^2 - 2xy - z^2) = n \quad (2)$$

**Proof**

If we replace the term $-2x^2y^2$ by $-4x^2y^2 + 2x^2y^2$ in equation **(1)**
we obtain,
$x^4 + y^4 + z^4 + 2x^2y^2 - 2y^2z^2 - 2z^2x^2 - 4x^2y^2 = n$;
$(x^2 + y^2 - z^2)^2 - (2xy)^2 = n$;
$(x^2 + y^2 - z^2 + 2xy)(x^2 + y^2 - z^2 - 2xy) = n$, which is equation (2)

Further,

$$\left[(x+y)^2 - z^2\right] \cdot \left[(x-y)^2 - z^2\right] = n;$$

$(x+y+z)(x+y-z)(x-y+z)(x-y-z) = n$; which is equation (3)

**Lemma 2**
*Let A, B, C be real numbers. Then, the three-variable linear system*
$$\begin{cases} x+y+z = A \\ x-y+z = B \\ x+y-z = C \end{cases}, \text{ has a unique solution; that}$$

solution being (x, y, z) = (x, y, z) = $\left(\dfrac{B+C}{2}, \dfrac{A-B}{2}, \dfrac{A-C}{2}\right)$

**Proof**

A routine calculation shows that the determinant of the matrix of the coefficients is equal to 4:

$D = \det \begin{bmatrix} 1 & 1 & 1 \\ 1 & -1 & 1 \\ 1 & 1 & -1 \end{bmatrix} = 4$, which is not zero. Therefore the given linear system has a unique solution: We apply Cramer's Rule. We have,

$D_x = \det \begin{bmatrix} A & 1 & 1 \\ B & -1 & 1 \\ C & 1 & -1 \end{bmatrix} = 2(B+C)$, $D_y = \det \begin{bmatrix} 1 & A & 1 \\ 1 & B & 1 \\ 1 & C & -1 \end{bmatrix} = 2(A-B)$, and

$D_z = \det \begin{bmatrix} 1 & 1 & A \\ 1 & -1 & B \\ 1 & 1 & C \end{bmatrix} = 2(A-C)$. Therefore the unique solution is given by,

$x = \dfrac{D_x}{D} = \dfrac{B+C}{2}$, $y = \dfrac{D_y}{D} = \dfrac{A-B}{2}$, $z = \dfrac{D_z}{D} = \dfrac{A-C}{2}$ ∎

**Lemma 3**

*Let A, B, C, D be real numbers. Consider the 3-variable linear system,*

$\begin{cases} x + y + z = A \\ x - y + z = B \\ x + y - z = C \\ x - y - z = D \end{cases}$

*(i) If the real numbers, A, B, C, D satisfy the condition,*
*A+D=B+C; then the given system has the unique solution*

$(x, y, z) = \left( \dfrac{B+C}{2}, \dfrac{A-B}{2}, \dfrac{A-C}{2} \right)$

(ii) If A+D is not equal to B+C; the above system has no solution

**Proof**

(i) A routine calculation shows that since A+D=B+C; the triple $\left( \dfrac{B+C}{2}, \dfrac{A-B}{2}, \dfrac{A-C}{2} \right)$

satisfies the last or fourth equation, which is $x - y - z = D$. By lemma 2, we also know that the same triple is the unique solution of the sub-system consisting of the first three equations. Since every solution of the system of four equations; is also a solution of any sub-system of equations; it follows that above triple is the unique solution of the system of four equations

(ii) Since any solution of the system oif four equations must be solution to the sub-system of the first three equations, it follows by Lemma 2, that the given system (of four equation) can have at most one solution, depending on whether the fourth equation is satisfied or not. Thus, if $A+D \neq B+C$; the fourth equation is not satisfied and so the system has no solution ∎

## 3. Theorems 1 and 2, and a solution to Contest problem CC24

**Theorem 1**

*Suppose that n is a positive integer which is exactly divisible by 8. In other words, n=8N, where N is an odd positive integer. Then the equation,*

$$x^2 + y^4 + z^4 - 2x^2y^2 - 2y^2z^2 - 2z^2x^2 = n \quad (1)$$

*has no integer solutions*

**Proof**

First observe that every integer solution to (1), must satisfy a necessary condition: two of x, y, z must be odd; the third even. To see why this is true, consider the possibilities:
Either all three x, y, z are odd, or all three are even, or two of them are even, the third odd; or two amond them are odd, the third even.
The first possibility is ruled out by considering (1) modulo 2:
The left-hand said would be odd, while n is even.
The second possibility is ruled out by considering (1) modulo 16:
The left-hand side would be congruent to zero modulo 16;
while $n = 8N \equiv 8 \pmod{16}$, since N is odd
The third possibility is also ruled out modulo 2; the left-handed side would be odd, while n is even.

We conclude that if (x, y, z) is a solution to (1); then two among x, y, z must be odd, the third even. Also, due to the symmetry of (1); if $(x_0, y_0, z_0)$ is a solution to (1); so are the other five permutations of this triple: the triples $(x_0, z_0, y_0), (y_0, x_0, z_0), (y_0, z_0, x_0), (z_0, x_0, y_0)$, and $(z_0, y_0, x_0)$, Therefore we may assume, (x and y odd, z even) in equation (1).

By Lemma 1, equation (1) is equivalent to (2);
$$(x^2 + y^2 + 2xy - z^2)(x^2 + y^2 - 2xy - z^2) = 8N \quad (3)$$

Since $x \equiv y \equiv 1 \pmod{2}$ and $z \equiv 0 \pmod{2}$. Both factors on the left-hand side of equation **(3)**, are even. But N is odd. Therefore (3) implies

Either $\begin{cases} x^2 + y^2 + 2xy - z^2 = 4d_1 \\ x^2 + y^2 - 2xy - z^2 = 2d_2 \end{cases}$ **(4a)**

Or alternatively, $\begin{cases} x^2 + y^2 + 2xy - z^2 = 2d_1 \\ x^2 + y^2 - 2xy - z^2 = 4d_2 \end{cases}$ **(4b)**

Where $d_1, d_2$ are integers such that,

$d_1 d_2 = N$

The integers $d_1, d_2$ (both positive or both negative) are odd integers, divisors of N; and whose product is N. If (4a) holds, then by subtracting the second equation from the first we obtain,

$4xy = 4d_1 + 2d_2;$
$2(xy - d_1) = d_2$
, which is impossible since $d_2$ is odd. The same type of contradiction arises in the case of (4b). The proof, by contradiction is complete. ∎

**Theorem 2**

Let n be a positive integer, and a, b, c, d integers such that $\begin{cases} a+d = b+c; \\ \text{and } abcd = n \end{cases}$. Then the rational triple $\left(\dfrac{b+c}{2}, \dfrac{a-b}{2}, \dfrac{a-c}{2}\right)$ is a solution to the equation $x^4 + y^4 + z^4 - 2x^2 y^2 - 2y^2 z^2 - 2z^2 x^2 = n$

And therefore (due to symmetry and even exponents), the triples $(\pm xo, \pm yo, \pm zo)$; are also rational solutions; where (xo,yo,zo) is a permulation of $\left(\dfrac{b+c}{2}, \dfrac{a-b}{2}, \dfrac{a-c}{2}\right)$; and with any of the eight possible combinations allowed.

**Proof**

Since $a+d = b+c$, it follows by Lemma 3(i) that the triple $(x, y, z) = \left(\dfrac{b+c}{2}, \dfrac{a-b}{2}, \dfrac{a-c}{2}\right)$ is the unique solution to system of equations,

$$\begin{cases} x+y+z = a \\ x-y+z = b \\ x+y-z = c \\ x-y-z = d \end{cases}$$

Therefore, by multiplying the four equations memberwisel it is clear that the same rational triple is a solution to the equation, $(x+y+z)(x-y+z)(x+y-z)(x-y-z) = abcd = n$; which is equivalent, by Lemma 1, to equation (1). Therefore the above rational triple is a solution to $x^4 + y^4 + z^4 - 2x^2 y^2 - 2y^2 z^2 - 2z^2 x^2 = n$ ∎

In the case of $n = 24 = 4 \cdot 3 \cdot 2 \cdot 1$; we have $4+1 = 3+2$. So, with $a = 4$, $d = 1$, $b = 3$, $c = 2$. We obtain the rational solution $(x, y, z) = \left(\dfrac{5}{2}, \dfrac{1}{2}, 1\right)$.

## 4. A remark

**Remark 1**

*Since equation (1) is symmetric with respect to the three variables x, y, z, and since all the exponents of the variables are even; it follows that if (xo, yo, zo) is a solution to equation (1); then so are the following forty eight (not necessarily distinct) triples (with (xo, yo, zo) being among them):*

$$(\pm|x_o|, \pm|y_o|, \pm|z_o|), (\pm|x_o|, \pm|z_o|, \pm|y_o|), (\pm|y_o|, \pm|x_o|, \pm|z_o|),$$
$$(\pm|y_o|, \pm|z_o|, \pm|x_o|), (\pm|z_o|, \pm|x_o|, \pm|y_o|), (\pm|z_o| \pm|y_o|, \pm|x_o|);$$

*With all eight sign combinations being possible in each of the six permutations.*

*Observe that if all three absolute values $|x_0|, |y_0|, |z_0|$ are distinct positive integers or rationals. Then the above 48 triples are distinct. If $|x_0|, |y_0|, |z_0|$ are positive with two of them being equal; the third different; then only twenty four of these triples are distinct. Likewise if one of $|x_0|, |y_0|, |z_0|$ is zero; the other two distinct positive integers. If one of $|x_0|, |y_0|, |z_0|$ is zero and the other two equal positive integers; then only twelve of the above 48 triples are distinct. If two among $|x_0|, |y_0|, |z_0|$ are zero; the third a positive integer, then only 6 are distinct. However, this last combination can occur only when n is the fourth power of a positive integer.*

### Theorem 3 and its proof

**Theorem 3**

*Let n be a positive integer. Then the 3-variable equation, $x^4 + y^4 + z^4 - 2x^2y^2 - 2y^2z^2 - 2z^2x^2 = n$; has an integer solution if, and only if, there exist positive integers a, b, c, d, which satisfy the three conditions*

$$\begin{cases} a+d = b+c, \ abcd = n, \\ \text{and } a \equiv b \equiv c \equiv d \pmod{2} \end{cases}$$

*(i.e. all four are even or odd)*

**Proof**

First, assume that a, b, c, d are positive integers satisfying the above three conditions. Then, the three rational numbers, $x_o = \dfrac{b+c}{2}$, $y_o = \dfrac{a-b}{2}$, $z_o = \dfrac{a-c}{2}$; are all integers since a, b, c are either all even; or are all three odd (by the third condition). Moreover, since a+d=b+c, it follows by Lemma 3(i), that the triple is $(x_0, y_0, z_0)$, the unique solution to the system,

$$\begin{cases} x+y+z=a \\ x-y+z=b \\ x+y-z=c \\ x-y-z=d \end{cases}; \text{ which in turn implies that } (x_0,\ y_0,\ z_0), \text{ is an integer solution of the equation,}$$

$(x+y+z)(x-y+z)(x+y-z)(x-y-z) = abcd = n$; and thus, by Lemma 1, a solution of equation **(1)**.

Next we prove the converse. Suppose that $(x_0,\ y_0,\ z_0)$, is an integer solution to equation (1). Then, by Lemma 1 it follows that $(x_0,\ y_0,\ z_0)$, is a solution of equation (3); and so,

$(x_o + y_o + z_o)(x_o - y_o + z_o)(x_o + y_o - z_o)(x_o - y_o - z_o) = n$   **(5)**

By (5) it follows that

$$\begin{cases} x_o + y_o + z_o = d_1 \\ x_o - y_o + z_o = d_2 \\ x_o + y_o - z_o = d_3 \\ x_o - y_o - z_o = d_4 \end{cases} \quad \textbf{(6a)}$$

And with the integers $d_1, d_2, d_3, d_4$ satisfying,

$d_1 d_2 d_3 d_4 = n$   **(6b)**

According to (6a), the triple $(x_0,\ y_0,\ z_0)$, is a solution to the system

$\{x+y+z=d_1,\ x-y+z=d_2,\ x+y-z=d_3,\ x-y-z=d_4\}$

Therefore, by Lemma 3 it follows that,

$\left( x_o = \dfrac{d_2 + d_3}{2},\ y_o = \dfrac{d_1 - d_2}{2},\ z_o = \dfrac{d_1 - d_3}{2} \right)$   **(7a)**

And also that, $d_1 + d_4 = d_2 + d_3$   **(7b)**

Since xo, yo, zo are integers; it follows from (7a); that all three integers $d_1, d_2, d_3$ must have the same parity; they must all be even or all odd; and so by (7b), $d_4$ must also have the same parity. Therefore,

$d_1 \equiv d_2 \equiv d_3 \equiv d_4 \pmod{2}$   **(8)**

Since n is a positive integer, then (6b) and (7b) we deduce that one of three possibilities must occur:
Possibility 1: All four $d_1, d_2, d_3, d_4$ are positive integers.
Possibility 2: All four are negative integers.

Possibility 3: One of $d_1, d_4$ is positive, the other negative. And likewise, one of $d_2, d_3$ is positive, the other negative.

If Possibility 1 holds, we are done since (6b), (7b), and (8) are the three conditions we seek to establish. If Possibility 2 holds; by setting $a = -d_1, b = -d_2, c = -d_3, d = -d_4$. Then (6b), (7b), and (8); imply $abcd=n$, $a+d=b+c$, and also $a \equiv b \equiv c \equiv d \pmod{2}$. We are done in this case. If Possibility 3 occurs. Then there are four possible combinations, all treated similarly. So, when $d_1 > 0$, $d_4 < 0$, $d_2 > 0$, $d_3 < 0$. Just set $a = d_1, c = -d_4, b = d_2, d = -d_3$; and we are done: (6b), (7b), and (8) imply, $abcd = n$, $a+d = b+c$, $a \equiv b \equiv c \equiv d \pmod{2}$ ∎

## 5. An Application of Theorem 3: Theorem 4

**Theorem 4**

*Let n be a positive integer, n=p, a prime; or n=4; or n=pq, a product of two distinct primes p and q. Then, the equation $x^4 + y^4 + z^4 - 2x^2y^2 - 2y^2z^2 - 2z^2x^2 = n$ has no integer solutions.*

**Proof**

According to Th.3, this equation will have an integer solution, if and only if there exist positive integers a,b,c,d satisfying the three conditions $\left\{\begin{array}{l} a+d = b+c, \ abcd = n, \\ \text{and } a \equiv b \equiv c \equiv d \pmod{2} \end{array}\right\}$

If n=p, a prime. Then $abcd=n=p$ implies that three of the positive integers a, b, c, d; must equal 1, the fourth 0. But then the additive condition $a+d = b+c$ cannot be satisfied, since one side will equal 2; the other side $p+1 \geq 3 > 2$. If $n = pq$; $abcd = pq$ implies that either one of a, b, c, d is pq, the other three 1; which renders a+d=b+c impossible, since one side will be $pq+1 \geq 2 \cdot 3 + 1 = 7 > 2 =$ other side. Or alternatively, one of $a,b,c,d$ is $p$, another is $q$; the remaining two, each being 1. And so $a+d = b+c$ implies either $1+p = 1+q$; or 2=p+q, both impossible since p, q are distinct primes.

Finally if, $n = 4$, a quick check shows that the first two conditions are satisfied only when one of a, d is 1, the other 2. And likewise, one of b, c is 1, the other 2; so $1+2 = 1+2$ and $1 \cdot 2 \cdot 1 \cdot 2 = 4$. However, the congruence condition obviously fails. We have shown that there exist no positive integers a, b, c, d satisfying all three conditions. Thus equation **(1)** has no integer solutions. ∎

## 6. Theorem 5 and it's proof

The following theorem provides a method for finding all the nonnegative integer solutions of equation **(1).**

**Theorem 5**

*Let n be a positive integer and consider the equation,*

$$x^4 + y^4 + z^4 - 2x^2y^2 - 2y^2z^2 - 2z^2x^2 = n$$

*Then, if $(x_0, y_0, z_0)$, is a nonnegative integer solution to this equation, then one of three sets of conditions must occur:*

*Either* $C_1 : \begin{cases} x_o = \dfrac{b+c}{2}, y_0 = \dfrac{a-b}{2}, z_0 = \dfrac{a-c}{2}, \\ \text{where a, b, c, d are positive integers} \\ \text{satisfying } a+d = b+c, \ abcd = n, \\ a \equiv b \equiv c \equiv d \pmod{2}; \text{ and} \\ a \geq b, \ a \geq c, \ a \geq d \end{cases}$

*Or* $C_2 : \begin{cases} x_0 = \dfrac{b-c}{2}, y_0 \dfrac{a-b}{2}, z_0 = \dfrac{a+c}{z}, \\ \text{where } a,b,c,d \text{ are positive integers} \\ \text{satisfying } a+d = b+c, \ abcd = n, \\ a \equiv b \equiv c \equiv d \pmod{2}; \text{ and} \\ a \geq b \geq c. \end{cases}$

*Or* $C_3 : \begin{cases} x_0 = \dfrac{c-b}{2}, y_0 = \dfrac{a+b}{2}, z_0 = \dfrac{a-c}{2} \\ \text{where } a,b,c,d \text{ are positive integers} \\ \text{satisfying } a+d = b+c, \ abcd = n, \\ a \equiv b \equiv c \equiv d \pmod{2}; \text{ and} \\ a \geq c \geq b \end{cases}$

Conversely, every nonnegative solution can be obtained from one of the three sets of conditions; $C_1, C_2,$ or $C_3$

**Proof**

First we prove the converse. Suppose that one of three sets of conditions is satisfied. If $C_1$ is satisfied, then it follows that $x_0, y_0, z_0$ are nonnegative integers, as $a \equiv b \equiv c \pmod 2$ and $a \geq b, a \geq c$ clearly imply. Also, since $a+d = b+c$, Lemma3(i) implies that $(x_0, y_0, z_0)$ is the unique solution to the system, $x+y+z = a, \ x-y+z = b, \ x+y-z = c, \ x-y-z = d.$ And, since $abcd = n,$ by Lemma 1, it follows that $(x_0, y_0, z_0)$ is a nonnegative integer solution of equation (1)

If $C_2$ is satisfied, a similar argument shows that $(x_0, y_0, z_0)$ is a nonnegative integer solution to (1).

Lemma 3(i) is applied with A=a, B=b, C=-c, and D=-d; which gives

$(x+y+z)(x-y+z)(x+y-z)\cdot(x-y-z) = a\cdot b\cdot(-c)\cdot(-d) = abcd = n$

If $C_2$ is satisfied, we apply Lemma 3(i) with A=a, B=-b, C=c and D=-d; $a\cdot(-b)\cdot c\cdot(-d) = abcd = n$.

Again $(x_0, y_0, z_0)$ is a nonnegative solution of equation (1)

Next, suppose that $(x_0, y_0, z_0)$ is a nonnegative solution of equation (1); and therefore by Lemma 1, of equation (3) as well: $(x_0+y_0+z_0)(x_0-y_0+z_0)(x_0+y_0-z_0)(x_0-y_0-z_0) = n$ **(9)**

Since $x_0 \geq 0, y_0 \geq 0, z_0 \geq 0$; we have $x_0 + y_0 + z_0 \geq 0$; in fact,

$$\begin{cases} a = x_0 + y_0 + z_0 > 0; \text{ with } x_0 \geq 0, y_0 \geq 0, z_0 \geq 0 \\ a = x_0 + y_0 + z_0 \geq 1; \\ \text{since not all nonnegative integers } x_0, y_0, z_0 \\ \text{can be zero; by (9) and the fact that} \\ \text{n is a positive integer.} \end{cases} \quad (10)$$

It follows from **(9)** and **(10)** that either all three integers $(x_0-y_0+z_0), (x_0+y_0-z_0), (x_0-y_0-z_0)$; are positive. Or otherwise, two of them must be negative; the third positive.

If all three $(x_0-y_0+z_0), (x_0+y_0-z_0), (x_0-y_0-z_0)$ are positive integers. Then,

$$\begin{cases} x_0 + y_0 + z_0 = a, x_0 - y_0 + z_0 = b, x_0 + y_0 - z_0 = c, x_0 - y_0 - z_0 = d \\ \text{with } a, b, c, d \text{ positive integers} \end{cases} \quad (10a)$$

Clearly, since $(x_0, y_0, z_0)$ are nonnegative; Also from **(10a)** implies $a \geq b$, and $a \geq c$. Also from **(10a)** and **(9)**, we have abcd=n. Moreover **(10a)** says that $(x_0, y_0, z_0)$ is a solution to the linear system,

$x+y+z=a, x-y+z=b, x+y-c=c, x-y-z=d$ Therefore by Lemma 3 it follows that

a+d=b+c; and that $x_0 = \dfrac{b+c}{2}, y_0 = \dfrac{a-b}{2}, z_0 = \dfrac{a-c}{2}$; which in turn implies that, since $x_0, y_0, z_0$ are integers; we must have $a \equiv b \equiv c \pmod{2}$; and by a+d=b+c; $a \equiv b \equiv c \equiv d \pmod{2}$; It is now clear that all the conditions in $C_1$ are satisfied.

If two of the factors $x_0 - y_0 + z_0, x_0 + y_0 - z_0, x_0 - y_0 - z_0$; are negative, the third positive; there are three possibilities. The possibility $x_0 - y_0 + z_0 < 0, x_0 + y_0 - z_0 < 0, x_0 - y_0 - z_0 > 0$; is ruled out; since

$x_0 - y_0 + z_0 < 0$ and $x_0 + y_0 - z_0 < 0$ imply $2x_0 < 0$ contrary to $x_0 \geq 0$.

Of the other two possibilities:

If $x_0 - y_0 + z_0 > 0, x_0 + y_0 - z_0 < 0, x_0 - y_0 - z_0 < 0$.

We then have,

$$\begin{cases} x_0 + y_0 + z_0 = a,\ x_0 - y_0 + z_0 = b,\ x_0 + y_0 - z_0 = -c,\ x_0 - y_0 - z_0 = -d \\ \text{with } a,b,c,d \text{ positive integers.} \end{cases} \quad \textbf{(10b)}$$

Since $y_0 \geq 0$, it follows from the first two equations in **(10b)** that $a \geq b$. And since $x_0 \geq 0$, it follows from the second and third equations in **(10b)** that $b \geq c$; so that $a \geq b \geq c$.

From **(10b)** and **(9)** we deduce that, $a \cdot b \cdot (-c) \cdot (-d) = n;\ abcd = n$

Next, apply Lemma 3 with $A=a,\ B=b,\ C=-c,\ D=-d$. It follows that

$$x_0 = \frac{b-c}{2},\ y_0 = \frac{a-b}{2},\ z_0 = \frac{a+c}{2},\ \text{and } A+D = B+C;\ a+d = b+c.$$

Also, since $x_0, y_0, z_0$ are integers; it follows that $a \equiv b \equiv c \pmod{2}$; and from $a+d = b+c$; we further have $a \equiv b \equiv c \equiv d \pmod{2}$.

Lastly, the remaining possibly is treated similarly as the previous one.

If $x_0 + y_0 - z_0 > 0,\ x_0 - y_0 + z_0 < 0,\ x_0 - y_0 - z_0 < 0$. We have,

$$\begin{cases} x_0 + y_0 + z_0 = a,\ x_0 - y_0 + z_0 = -b,\ x_0 + y_0 - z_0 = c,\ x_0 - y_0 - z_0 = -d \\ \text{where } a,b,c,d \text{ are positive integers} \end{cases}$$

It follows that (since $z_0 \geq 0$) $a \geq c$;

and that $c \geq b$ (since $z_0 \geq 0$);

And so $a \geq c \geq b$

We also have $a \cdot (-b) \cdot c \cdot (-d) = abcd = n$

After that, we apply Lemma3 with,

$A=a,\ B=-b,\ C=c,\ D=-d.$

We omit the details. ∎

### 7. An application of Theorem 5: Theorem 6

Consider the case $n = p^2$, $p$ an odd prime. We can use Theorem 5 to determine all the nonnegative integer solutions of equation **(1)**. Observe that if four positive integers $a,b,c,d$ satisfy the conditions,

$abdc = p^2$ and $a+d = b+d$

Then since $p$ is prime; one of $a$, $d$ must equal $p$, the other 1; and one of $b$, $c$ must equal $p$, the other 1. And since $p$ is odd, all four $a, d, b, c$ are odd, so the congruence condition $a \equiv b \equiv c \equiv d \pmod{2}$ is satisfied.

Also note that the possibility that one of $a, b, c, d$ is $p^2$; the other three equal to 1; is obviously ruled out by the condition $a+d=b+c$.

We know from Theorem 5, that every nonnegative solution can be obtained from one of the sets of conditions $C_1, C_2,$ or $C_3$. Also note that in all three $C_1, C_2, C_3$; $a = \max\{a,b,c\}$. This clearly implies that always, a=p. Therefore:

If $C_1$ is satisfied, then a=p, d=1; and either b=p, c=1; or b=1 and c=p. We obtain two nonnegative solutions: $\left(\frac{p+1}{2}, \frac{p-1}{2}, 0\right)$ $\left(\frac{p+1}{2}, 0, \frac{p-1}{2}\right)$. If $C_2$ is satisfied, we have $a \geq b \geq c$. Therefore, $a = p, d = 1, b = p, c = 1$. And therefor according to $C_2$, we obtain the nonnegative solution $\left(\frac{p-1}{2}, 0, \frac{p+1}{2}\right)$.

If $C_3$ is satisfied, we have $a \geq c \geq b$. Thus, a=p, d=a, b=1, c=p. We obtain the nonnegative solution $\left(\frac{p-1}{2}, \frac{p+1}{2}, 0\right)$.

From these four nonnegative integer solutions, all integer solutions of equation (1) can be obtained from permutations and combinations of signs. We have the following.

**Theorem 6**

Let $n = p^2$, p an odd prime. Then, the 3-variable equation

$x^4 + y^4 + z^4 - 2x^2y^2 - 2y^2z^2 - 2z^2x^2 = n,$ has
exactly twenty four integer solutions:
$(\pm e, \pm f, 0), (\pm e, 0, \pm f), (0, \pm e, \pm f),$
$(0, \pm f, \pm e), (\pm f, \pm e, 0), (\pm f, 0, \pm e);$
where $e = \frac{p+1}{2}, f = \frac{p-1}{2}.$

*And in each triple, all four sign combinations can hold.*

8. **When n is a perfect square: $n = k^2$**

In this section, we determine all the integer solutions, with xyz=0, in the case $n = k^2$, k a positive integer. We also determine, when $n = k^2$, all the rational solutions of equation (1), with xyz=0. To accomplish our goal we need three results.

**Result 1**

*Let k be a positive integer $k \geq 2$ Consider the 2-variable equation, $u^2 - v^2 = k$*
*(i) If $k \equiv 2 \pmod 4$, then this equation has no positive integer solutions.*
*(ii) If $k \equiv 0, 1, or 3 \pmod 4$. Then this equation has positive integer solutions. The solution set, in positive*

*integers u and v can be parametrically described as follows:*

$$\begin{cases} u = \dfrac{d_1 + d_2}{2}, \ v = \dfrac{d_1 - d_2}{2}; \text{ where } d_1, d_2 \text{ are} \\ \text{positive integers such that } d_1 > d_2, \ d_1 \equiv d_2 \pmod{2}, \\ \text{and } d_1 d_2 = k \end{cases}$$

*Each positive integer solution pair (u,v) can be obtained once in this way: if $D_1, D_2$ are positive integers such that $D_1 > D_2$, $D_1 \equiv D_2 \pmod 2$, $D_1 D_2 = k$; and with $d_1 \neq D_1$ or $d_2 \neq D_2$. Then the integer solution pair (U,V) is distinct from the solution pair (u,v); where $U = \dfrac{D_1 + D_2}{2}$ and $V = \dfrac{D_1 - D_2}{2}$.*

**Proof**

(i) An integer square can only be congruent to 0 or 1 modulo 4; according to whether that integer is even or odd. Therefore,

$u^2 - v^2 \equiv 1-1, 1-0, 0-1, \text{ or } 0-0 \pmod 4;$

$u^2 - v^2 \equiv 0, 1, -1; \text{ or } 0 \pmod 4;$

$u^2 - v^2 \equiv 0, 1, \text{ or } 3 \pmod 4.$

Therefore $u^2 - v^2$ cannot be congruent to 2mod4; and so there are no integer solutions.

(ii) If $k \equiv 0 \pmod 4$ then k can be written as a product of two positive even integers:

$k = u^2 - v^2 = (u-v)(u+v);$ and so we must have,

$u + v = d_1, u - v = d_2 (u > v, \text{ since } u, v, \in \mathbb{Z}^+, \text{ and } k \in \mathbb{Z}^+)$ Where $d_1, d_2$ are even positive integers with $d_1 > d_2$, and $d_1 d_2 = k$. From which we obtain, $u = \dfrac{d_1 + d_2}{2}, v = \dfrac{d_1 - d_2}{2}$

If $k \equiv 1$ or $3 \pmod 4$; then k is an odd integer; $k = d_1 d_2$, with $d_1 \equiv d_2 \equiv 1 \pmod 2$, $d_1 > d_2 \geq 1$. Again,

$u = \dfrac{d_1 + d_2}{2}, v = \dfrac{d_1 - d_2}{2}$. That each solution is obtained in this way only once is fairly obvious; we omit the details. ∎

The following result states the well-known parametric formulas that describe all the positive integer solutions to the Pythagorean equation. This result can be found in number theory books, in wikipedia on the world wide web, and other internet sources. For a wealth of information regarding Pythagorean triples and triangles, an excellent source can be found in reference [2]. Also, in the same reference, a detailed derivation of the parametric formulas stated below.

**Result 2**

*All the positive integer solutions to the 3-variable equation $u^2 + v^2 = w^2$ can be parametrically described as follows: $u = d(k_1^2 - k_2^2), v = d \cdot (2k_1k_2), w = d(k_1^2 + k_2^2)$; or alternatively,*
*$u = d(2k_1k_2), v = d(k_1^2 - k_2^2), w = d(k_1^2 + k_2^2)$. Where $d, k_1, k_2$ are positive integers such that*
*$k_1 > k_2, k_1 + k_2 \equiv 1 \pmod{2}$ (i.e. one of $k_1, k_2$ is odd, the other even); and $\gcd(k_1, k_2) = 1$ (i.e. $k_1$ and $k_2$ are relatively prime.*

In a 2006 paper published in the journal
*Mathematics and Computer Education*         (see reference [3]),

By this author; the reader can find parametric formulas describing the entire set of positive integer solutions to the 3-variable equation $u^2 + kv^2 = w^2$; where $k \geq 2$, is a fixed positive integer. In the same article, the reader can find a step-by-step derivation of all the positive integer solutions to this equation. We state this result.

**Result 3**

*Let k be a positive integer, $k \geq 2$. The entire set of positive integer solutions of the 3-variable equation $u^2 + k \cdot v^2 = w^2$, can be parametrically described as follows:*

$$u = \frac{d \cdot |k_1 \cdot m_1^2 - k_2 \cdot m_2^2|}{2}, v = dm_1m_2, w = \frac{d \cdot (k_1 \cdot m_1^2 + k_2 \cdot m_2^2)}{2}$$

*Where $d, k_1, k_2, m_1, m_2$ are positive integers such that*
*$k_1 \cdot k_2 = k$, $\gcd(m_1, m_2) = 1$, (i.e. $m_1$ and $m_2$ are relatively prime) and $dk_1m_1^2 \equiv dk_2m_2^2 \pmod{2}$*

Consider the case where two among *x, y, z* in equation **(1)**, are zero; say *y=z=0*. Then,
$x^4 = k^2 \Leftrightarrow |x| = \sqrt{k}$; which shows that if k is a perfect square; $k = m^2$; with $m \in \mathbb{Z}^+$,
then $|x| = m$, and so there are integer solutions; otherwise there are no integer solutions. What about rational solutions? Recall that the square root of a positive integer is either a positive integer; or otherwise an irrational number. This statement is equivalent to the following statement: A positive integer is equal to the square of a rational number; if, and only if, that same integer is a perfect or integer square. This holds true for any natural number exponent, and not just squares. This statement can be fairly easily proven. The interested reader may refer to reference **[2]**.

**Theorem 7**

*Let $n = k^2$, k a positive integer, and consider the equation*
*$x^4 + y^4 + z^4 - 2x^2y^2 - 2x^2y^2 - 2y^2z^2 - 2z^2x^2 = n$. Let $S_1$ and $S_2$ be the sets,*
*$S_1 = \{(a,b,c) | (a,b,c) \text{ is an integer solution (1) with two among a,b,c being zero}\}$; and*
*$S_2 = \{(a,b,c) | (a,b,c) \text{ is a rational solution of equation (1), with two among } a,b,c \text{ being zero}\}$*

*Then,*

*(i) If K is not a perfect or integer square, $S_1 = S_2 = \emptyset$; equation (1) has no such integer or rational solutions.*

*(ii) If $k = m^2$, m a positive integer; then,*
*$S_1 = S_2 = \{(\pm m, 0, 0), (0, \pm m, 0), (0, 0, \pm m)\}$ There are six integer and six identical rational solutions in this case.*

Next, we consider the case where exactly one of the three variables is zero; the other two nonzero. If we take z=0, and $xy \neq 0$. Then $(x^2 - y^2)^2 = k^2$; $|x^2 - y^2| = k$; with $k \geq 2$. And because of the symmetry of equation (1), to find all such integer solutions; it suffices to find all the positive integer solutions of the equation, $x^2 - y^2 = k$.

But Result 1, provides a complete answer to this question. Accordingly, we can state Theorem 8.

**Theorem 8**

*Let $n = k^2$, k an integer with $k \geq 2$. Consider the equation, $x^4 + y^4 + z^4 - 2x^2y^2 - 2y^2z^2 - 2z^2x^2 = n$*
*(1)*

*And let S be the set,*
*$S = \{(a,b,c) | (a,b,c)$ is an integer solution of equation (1), with exactly one of a,b,c being zero$\}$*
*(i) If $k \equiv 2 \pmod 4$, $S = \emptyset$; there are no such solutions.*
*(ii) If $k \equiv 0, 1$, or $3 \pmod 4$; then the set S is a finite nonempty set which can be described/determined in the following way:*
*For each pair of divisors $d_1$ and $d_2$ of k; such that $d_1 \equiv d_2 \pmod 2, 1 \leq d_2 < d_1$ and $d_1 d_2 = k$; define the positive integers $e = \dfrac{d_1 + d_2}{2}$ and $f = \dfrac{d_1 - d_2}{2}$ (The integers $d_1$ and $d_2$ satisfy the conditions of Result 1).*
*From a specific such pair $(d_1, d_2)$ of divisors $d_1$ and $d_2$ of k; there are exactly twenty four distinct integer solutions (with exactly one of x,y,z being zero) to equation (1):*
*$(\pm e, \pm f, 0), (\pm e, 0, \pm f), (\pm f, \pm e, 0), (\pm f, 0, \pm e), (0, \pm e, \pm f), (0, \pm f, \pm e)$. The number of integer solutions in the set S is 24N; where N is the number of the divisor pairs $(d_1, d_2)$ that satisfy the conditions, $d_1 \equiv d_2 \pmod 2$, $1 \leq d_2 < d_1$, and $d_1 d_2 = k$ (it is clear from Result 1, that the twenty four such solutions generated by another divisor pair $(D_1, D_2)$ will be distinct from the twenty four solutions generated by the pair $(d_1, d_2)$.*

Now, let us look at the case of rational solutions (with exactly one of the variables being zero) of equation (1). If we set z=0, we have $(x^2 - y^2)^2 = k^2$; and due to symmetry in (1), it suffices to find the positive rational solutions of $x^2 - y^2 = k$. By putting, $x = \dfrac{w}{v}, y = \dfrac{u}{v}$; where w, v, and u are positive integers.

There we obtain,

$u^2 + k \cdot v^2 = w^2$, which is the equation of Result 3.

Accordingly we get,

$$x = \frac{w}{v} = \frac{1}{2}\left(k_1 \cdot \frac{m_1}{m_2} + k_2 \cdot \frac{m_2}{m_1}\right); \text{ and}$$

$$y = \frac{u}{v} = \frac{1}{2} \cdot \left|\frac{k_1 m_1}{m_2} - \frac{k_2 m_2}{m_1}\right|; \text{ where } k_1, k_2 \text{ are positive integers such that } k_1 k_2 = k. \text{ And } m_1, m_2 \text{ are}$$

relatively prime positive integers.

**Theorem 9**

Let $n = k^2$, k an integer with $k \geq 2$. Consider the equation, $x^4 + y^4 + z^4 - 2x^2 y^2 - 2y^2 z^2 - 2z^2 x^2 = n$
*(1)*
*Let*

$S = \{(a,b,c) | (a,b,c) \text{ is a rational solution to equation (1), with exactly one of a, b, c being zero}\}$

*The set S is an infinite set which is the union of six classes of solutions; which can be parametrically be described as follows:*

*Let $k_1, k_2$ be positive integers such that $k_1 k_2 = k$ And $m_1, m_2$ (arbitrary) relatively prime positive integers.*

*Also, $e = \dfrac{1}{2} \cdot \left(\dfrac{k_1 m_1}{m_2} + \dfrac{k_2 m_2}{m_1}\right)$ (note that $e > 0$)*

$f = \dfrac{1}{2} \cdot \left|\dfrac{k_1 m_1}{m_2} - \dfrac{k_2 m_2}{m_1}\right|$ *And also with $k_1 m_1^2 \neq k_2 m_2^2$ (so that $f > 0$)*

*Class $C_1$: The set of rational triples of the form*

$(\pm e, \pm f, 0)$

*Class $C_2$: The set of rational triples of the form*

$(\pm e, 0, \pm f)$

*Class $C_3$: The set of rational triples of the form*

$(\pm f, \pm e, 0)$

*Class $C_4$: The set of rational triples of the form*

$(\pm f, 0, \pm e)$

*Class $C_5$: The set of rational triples of the form*

$(0, \pm e, \pm f)$

*Class $C_6$: The set of rational triples of the form*

$(0, \pm f, \pm e)$

*Then, $S = C_1 U C_2 U C_3 U C_4 U C_5 U C_6$*

*(With all sign combinations being possible)*

We finish this section with the case n=1. The case of integer solutions with xyz=0 is immediate: $(\pm 1, 0, 0), (0, \pm 1, 0), (0, 0, \pm 1)$, are the only integer solutions with two of the variables being zero. The same is true for such rational solutions. When only one of the variables is zero; there are no integer solitions since (as it is easily verified), the equation $x^2 - y^2 = 1$ has no positive integer solutions. The same equation though, has positive rational solutions. By setting, $x = \dfrac{w}{v}$, $y = \dfrac{u}{v}$; where w, v, u are positive integers. We obtain $w^2 = v^2 + u^2$. And from Result 2 we further get,

$$x = \frac{w}{v} = \frac{k_1^2 + k_2^2}{2k_1 k_2}, \quad y = \frac{u}{v} = \frac{k_1^2 - k_2^2}{2k_1 k_2}; \quad \text{or alternatively,} \quad x = \frac{k_1^2 + k_2^2}{k_1^2 - k_2^2}, \quad y = \frac{2k_1 k_2}{k_1^2 - k_2^2}.$$

Where $k_1, k_2$ are relatively prime positive integers such that $1 \leq k_2 < k_1$ and $k_1 + k_2 \equiv (\bmod\, 2)$ We have the following theorem.

**Theorem 10**

*Consider the equation, $x^4 + y^4 + z^4 - 2x^2 y^2 - 2y^2 z^2 - 2z^2 x^2 = 1$.*

*(a) This equation has six integer solutions with two of the variables being zero:*
$(\pm 1, 0, 0), (0, \pm 1, 0),$ and $(0, 0, \pm 1)$. *These are also the only such rational solutions.*

*(b) The above equation has no integer solutions with only one of the variables being zero.*

*(c) The set S of rational solutions with only one of the variables being zero; is the union of twelve classes of solutions which can be described as follows:*

*Let $k_1, k_2$ be positive integers such that $1 \leq k_2 < k_1, k_1 + k_2 \equiv 1 (\bmod\, 2)$; and $k_1, k_2$ are relatively prime.*

*Let $e = \dfrac{k_1^2 + k_2^2}{2k_1 k_2}, f = \dfrac{k_1^2 - k_2^2}{2k_1 k_2}, g = \dfrac{k_1^2 + k_2^2}{k_1^2 - k_2^2},$ and $h = \dfrac{2k_1 k_2}{k_1^2 - k_2^2}.$*

*(Clearly all four; e, f, g, and h; are positive)*

*Class $C_1$: The set of rational triples of the form*

$(\pm e, \pm f, 0)$

*Class $C_2$: The set of rational triples of the form*

$(\pm e, 0, \pm f)$

*Class $C_3$: The set of rational triples of the form*

$(\pm f, \pm e, 0)$

*Class $C_4$: The set of rational triples of the form*

$(\pm f, 0, \pm e)$

*Class $C_5$: The set of rational triples of the form*

$(0, \pm e, \pm f)$

*Class $C_6$: The set of rational triples of the form*

$(0, \pm f, \pm e)$

*For Classes $C_7$ through $C_{12}$; just replace e by g ; and f by h; and repeat the process.*
*(with all sign combinations being possible)*

**Remark**

*Note that Th. 9, there is the condition $k_1 m_1^2 \neq k_2 m_2^2$, to ensure that f is not zero; and so a positive integer. One may ask the question: when is $k_1 m_1^2 = k_2 m_2^2$; or equivalently, $\dfrac{k_1}{k_2} = \dfrac{m_2^2}{m_1^2}$? Since $m_1$ and $m_2$ are relatively prime; that would imply $k_1 = \delta \cdot m_2^2$ and $k_2 = \delta \cdot m_1^2$, for some positive integer $\delta$. And since $k_1 k_2 = k$; this in turn would imply $(\delta m_1 m_2)^2 = k$; which can only occur when k is a perfect square, $k = N^2, N$ a positive integer. And so in such a case, $\delta m_1 m_2 = N$; so the product of the two positive integer parameters $m_1$ and $m_2$, would be a divisor of N.*

9. **The equation $x^2 + y^4 + z^4 - 2xy^2 - 2xz^2 - 2y^2 z^2 = n$**

If we replace $x^2$ by $x$ in equation **(1),** we obtain an equation which is quadratic in $x$:

$x^2 - 2(y^2 + z^2)x + y^4 + z^4 - 2y^2 z^2 = n$; or equivalently

$x^2 - 2(y^2 + z^2)x + (y^2 - z^2)^2 - n = 0$     **(11)**

Equation **(11)** is symmetric only with respect to the variables y and z. Applying the quadratic formula yields, $x = y^2 + z^2 \pm \sqrt{(2yz)^2 + n}$     **(12)**

10. **The equation $x^2 + y^4 + z^4 - 2xy^2 - 2xz^2 - 2y^2 z^2 = n$:**
    **Integer Solutions**

It is clear from **(12)** that equation **(11)** will have integer solutions if and only if the integer $(2yz)^2 + n$ is an integer square:

$$\begin{cases} (2yz)^2 + n = T^2 \\ \text{where T is a nonnegative integer. But since } n \geq 1 \\ \text{and } (2yz)^2 \geq 0;\text{ it is clear that T must} \\ \text{be a positive integer} \end{cases} \quad (13)$$

The following analysis will quickly lead us in formulating Theorem 11.

1) By inspection, by looking at equation (11); we see that if the positive integer n is not a perfect square; the (11) cannot have any integer solutions with $xyz = 0$. Indeed, if $x = 0$ in (11); then $n = (y^2 - z^2)^2$. And if $y = 0$ or $z = 0$; then $n = x^2$.

2) If n is a perfect square; $n = k^2$; k a positive integer.

   a) Solutions with $x = 0$ require that, $(y^2 - z^2)^2 = k^2$; and because of the symmetry with respect to y and z; it suffices to determine all the integer solutions of $y^2 - z^2 = k$
   
   If $k = m^2$, m a positive integer. Then the solutions with $z = 0$; have $y = \pm m$.
   
   If k is not a perfect square; then there are no solutions with $yz = 0$.
   
   If $k \equiv 2 \pmod 4$, there are no integer solutions.
   
   If $k \equiv 0, 1,$ or $3 \pmod 4$; the integer solutions of $y^2 - z^2 = k$; with $yz \neq 0$; are determined by the positive integer solutions of $y^2 - z^2 = k$; just apply Result 1.
   
   b) Solutions with $x \neq 0$ and $yz = 0$. From **(13)**, we get $T = k$. And so for $y = z = 0$, we have **(from(12))** $x = \pm k$; and for $y \neq 0$, $z = 0$ we obtain $x = y^2 \pm k$; where y can be integer; with the exception that if $k = m^2$; $y \neq \pm m$, if the minus sign holds in $x = y^2 \pm k$; since $x \neq 0$.

3) Solutions with $yz \neq 0$

   As we have seen in #1 above; unless n is a perfect square; $x \neq 0$ as well. If n is a perfect square, a particular combination of y and z may result in $x = 0$.

Regardless, solutions with $yz \neq 0$; can be completely determined from the positive integer solutions of equation (13): $T^2 - (2yz)^2 = n$ According to Result 1, if $n \equiv 2 \pmod 4$, there are no positive integer solutions. Otherwise, we must have $T = \dfrac{d_1 + d_2}{2}$, $2yz = \dfrac{d_1 - d_2}{2}$; where $d_1$ and $d_2$ are positive integers such that, $1 \leq d_2 < d_1$, $d_1 \equiv d_2 \pmod 2$, $d_1 d_2 = n$.

The second equation above states that, $4yz = d_1 - d_2$; which implies $d_1 \equiv d_2 \pmod 4$ If $n \equiv 0 \pmod 4$, the combination $d_1 = \dfrac{n}{2}, d_2 = 2$; ensures that this equation is always satisfied for some divisors $d_1$ and $d_2$ of n. Likewise, when $n \equiv 1 \pmod 4$; one such combination is $d_1 = n$ and $d_2 = 1$; so that
$d_1 - d_2 = n - 1 \equiv 0 \pmod 4$ (When n is prime, this is the only such combination) However, for
$n \equiv 3 \pmod 4$, no such $d_1$ and $d_2$ can exist. Indeed, $d_1 - d_2 \equiv 0 \pmod 4$ and $d_1 d_2 = n \equiv 3 \pmod 4$ And
so, $d_1 \equiv d_2 \pmod 4$ which implies, $d_1 d_2 \equiv d_1 \cdot d_1 \equiv (d_1)^2 \equiv n \equiv 3 \pmod 4$; which is impossible since in
fact $d_1^2 \equiv 1 \pmod 4$, by virtue of the fact that $d_1$ is odd. Thus, equation (11) has no integer solutions with
$yz \neq 0$ if $n \equiv 3 \pmod 4$.

**Theorem 11**

Let n be a positive integer. Consider the equation $x^2 + y^4 + z^4 - 2xy^2 - 2xz^2 - 2y^2z^2 = n$

*1)* If n is not a perfect square, then the above equation has no integer solutions with $xyz = 0$; that is, with
at least one of x, y, or z being zero.

*2)* Suppose that $n = k^2$, k a positive integer.

*a)* If k is not a perfect square. Then the above equation has no integer solutions with x=0; and with one of
y, z also being zero.
*b)* If $k = m^2$, m a positive integer; then the integer solutions with x=0; and one of y. z also zero are:
$(0, \pm m, 0)$, $(0, 0, \pm m)$; four such solutions in total.
*c)* If $k \equiv 2 \pmod 4$, there are no integer solutions with x=0.
*d)* The integer solutions with $x \neq 0$ and yz=0 can be described as follows:
$x = t^2 \pm k$, $y = t$, $z = 0$; being one family of such solutions.
And $x = t^2 \pm k$, $y = 0$, $z = t$; being the other. Where t can be any integer.
And with one exception: If $k = m^2$, $m \in \mathbb{Z}^+$; and the minus sign holds; then
$t \neq \pm m; |t| \neq m$; so that $x \neq 0$

*3)* The integer solutions with $yz \neq 0$ can be described as follows:
If $n \equiv 2$ or $3 \pmod 4$, then there are no integer solutions. If $n \equiv 0$ or $1 \pmod 4$; then all such integer
solutions can be described in the following way. For every pair of divisors $d_1$ and $d_2$ such that

$d_1 d_2 = n$, $1 \leq d_2 < d_1$, $d_1 \equiv d_2 \pmod 4$ (when $n \equiv 0 \pmod 4$); $d_1 = \dfrac{n}{2}$, $d_2 = 2$ is such a combination.

For $n \equiv 1 \pmod 4$; $d_1 = n$, $d_2 = 1$ is such a combination); define $e = \dfrac{d_1 - d_2}{4}$ and $f = \dfrac{d_1 + d_2}{2}$ Then all the integer solutions determined by this pair of divisors are given by $x = \pm \left[ \rho^2 + \left( \dfrac{e}{\rho} \right)^2 \pm f \right]$, $y = \pm \rho$, $z = \pm \dfrac{e}{\rho}$; *where $\rho$ is a positive integer divisor of the positive integer e; and with all sign combinations being possible.*

The following theorem is the immediate consequence of parts 1 and 3 of Theorem 11.

**Theorem 12**

*Let n be a positive integer with $n \equiv 2$ or $3 \pmod 4$. Then the 3-variable equation,*
$x^2 + y^4 + z^4 - 2xy^2 - 2xz^2 - 2y^2z^2 = n$ *has no integer solutions.*

**11. The equation $x^2 + y^4 + z^4 - 2xy^2 - 2xz^2 - 2y^2z^2 = n$ : Rational solutions**

As we have already seen this equation is equivalent to **(11),** whose solutions are given by **(12):**
$x = y^2 + z^2 \pm \sqrt{(2yz)^2 + n}$

First, consider rational solutions with $yz = 0$. Clearly, such solutions exist precisely when n is a rational square; and since n is a positive integer, n is a rational only if it is a perfect or integer square. So, if $n = k^2$, k a positive integer; all the rational solutions with yz=0, are given by

$(\pm k, 0, 0)$, $(r^2 \pm k, \pm r, 0)$, $(r^2 \pm k, 0, \pm r)$; where r is arbitrary positive rational; and will all the sign combinations being possible.

Consider rational solutions with $yz \neq 0$. If we put $|2yz| = \dfrac{u}{v}$, where u and v are positive integers. Then,

$\sqrt{(2yz)^2 + n} = \sqrt{\dfrac{u^2 + nv^2}{v^2}} = \dfrac{\sqrt{u^2 + nv^2}}{v}$, which will be rational if and only if, $u^2 + nv^2 = w^2$, for some positive integer w. By result3, we must have

$$\begin{cases} u = \dfrac{d\left|n_1 t_1^2 - n_2 t_2^2\right|}{2}, \ v = dt_1 t_2, \text{ and} \\ w = \dfrac{d \cdot \left(n_1 t_1^2 + n_2 t_2^2\right)}{2}; \text{ where } n_1, n_2, t_1, t_2 \\ \text{are positive integers such that } n_1 n_2 = n, \\ \gcd(t_1, t_2) = 1, \text{ and } dn_1 t_1^2 \equiv dn_2 t_2^2 \pmod 2 \\ \text{And with } n_1 t_1^2 \neq n_2 t_2^2 \text{ (so that } u > 0\text{)} \end{cases} \quad (14)$$

If we take $|y| = r$, a positive rational. Then from $|2yz| = \dfrac{u}{v}$ we get $|z| = \dfrac{u}{2rv}$; and thus by **(14)** we have,

$|z| = \dfrac{1}{4r} \cdot \left| \dfrac{n_1 t_1}{t_2} - \dfrac{n_2 t_2}{t_1} \right|$. Therefore from **(12)** we further get, $x = y^2 + z^2 \pm \dfrac{w}{v}$; and by **(14)**

$x = r^2 + \dfrac{1}{16r^2} \cdot \left( \dfrac{n t_1}{t_2} - \dfrac{n_2 t_2}{t_1} \right)^2 \pm \dfrac{1}{2} \left( \dfrac{n_1 t_1}{t_2} + \dfrac{n_2 t_2}{t_1} \right)$

**Theorem 13**

*Let n be a positive integer. Consider the equation $x^2 + y^4 + z^4 - 2xy^2 - 2xz^2 - 2y^2 z^2 = n$*

*(a) If n is not a perfect square, this equation has no rational solutions with yz=0.*

*(b) If $n = k^2$, k a positive integer, then all the rational solutions with yz=0 are triples of the form:*
$(\pm k, 0, 0), \; (r^2 \pm k, \pm r, 0), \; \text{and} \; (r^2 \pm k, 0, \pm r);$ *where r is an arbitrary positive rational; and with any sign combination being possible.*

*(c) All the rational solutions with $yz \neq 0$, can be parametrically described as follows:*

$x = r^2 + \dfrac{1}{16r^2} \left( \dfrac{n t_1}{t_2} - \dfrac{n_2 t_2}{t_1} \right)^2 \pm \dfrac{1}{2} \left( \dfrac{n_1 t_1}{t_2} + \dfrac{n_2 t_2}{t_1} \right)$

$y = \pm r, \; \text{and} \; z = \pm \dfrac{1}{4r} \left| \dfrac{n_1 t_1}{t_2} - \dfrac{n_2 t_2}{t_1} \right|$

*Where r is an arbitrary positive rational, any sign combination being possible; and with $n_1, t_1, n_2, t_2$ being positive integers such that $n_1 n_2 = n$ and $\gcd(t_1, t_2) = 1$ And with $n_1 t_1^2 \neq n_2 t_2^2$.*